\documentclass[reqno]{amsart}
\usepackage{amssymb}
\usepackage{hyperref}

\begin{document}
\title[ Riemannian manifolds admitting a projective semi-symmetric ...]
{Riemannian manifolds admitting a projective semi-symmetric connection}

\author[S. K. Chaubey$^{1}$, S. K. Yadav$^{2}$ and Pankaj$^{3}$]
{S. K. Chaubey$^{1}$, S. K. Yadav$^{2}$ and Pankaj$^{3}$} 

\address
 {$^{1}$
 Section of Mathematics, Department of Information Technology, Shinas College of Technology, Shinas, P.O. Box 77
Postal Code 324,  Oman.}
\email {sk22$_{-}$math@yahoo.co.in}

\address
{$^{2}$
Department of Mathematics, Poornima college of Engineering, ISI-6, RIICO, Institutional Area, Sitapura, Jaipur-302022, Rajasthan, India.}

\email {prof$_{-}$sky$16$@yahoo.com}

\address
 {$^{3}$
 Section of Mathematics, Department of Information Technology, Higher College of Technology, Muscat,  Oman.}
\email {pankaj.fellow@gmail.com}



\begin{abstract}
The aim of the present paper is to study the properties of Riemannian manifolds equipped with a projective semi-symmetric connection.
\end{abstract}

\maketitle
\numberwithin{equation}{section}
\newtheorem{theorem}{Theorem}[section]
\newtheorem{lemma}[theorem]{Lemma}
\newtheorem{proposition}[theorem]{Proposition}
\newtheorem{corollary}[theorem]{Corollary}
\newtheorem*{remark}{Remark}
\newtheorem{Agreement}[theorem]{Agreement}
\newtheorem{definition}[theorem]{Definition}

\section{Introduction} 
At the foundation of Riemannian geometry there are three ideas. The first of these is the realization of the fact that a non-Euclidean geometry exists (N. I. Lobachevskii). The second is the concept of the interior geometry of surfaces by C. F. Gauss and third is the concept of an $n-$dimensional space by B. Riemann. Our present paper belongs to the study of third case. The idea of Riemannian geometry played an important role in the formulation of the general theory of relativity by A. Einstein. A non-flat $n-$dimensional Riemannian manifold $(M_n, g)$, $(n>2)$, is said to be quasi Einstein manifold if its Ricci tensor $S$ of type $(0, 2)$ is not identically zero and satisfies the tensorial expression
\begin{equation}
\label{1.1}
S(X,Y)=ag(X,Y)+b\pi(X)\pi(Y),{\hspace{15.5pt}} X, Y \in{TM}
\end{equation}
for smooth functions $a$ and $b (\neq{0})$, where $\pi$ is a non-zero $1-$form associated with the  Riemannian metric $g$         
and the associated unit vector field $\xi$ \cite{chacki1}. The unit vector field $\xi$ is called the generator  and the $1-$form $\pi$ is called the associated $1-$form of the manifold. It is observed that a collection of non-interacting pressure less perfect fluid of general relativity is a four dimensional semi-Riemannian quasi Einstein manifold whose associated scalars are $\frac{r}{2}$ and $\kappa \rho$, where $\kappa$ is the gravitational constant, $\rho$ and $r$ are the energy density and scalar curvature, the generator of the manifold being the time like velocity vector field of the perfect fluid. If the generator of a quasi Einstein manifold is parallel vector field, then the manifold is locally a product manifold of  one-dimensional distribution $U$ and $(n-1)$ dimensional distribution $U^{\perp}$, where $U^{\perp}$ is involutive and integrable \cite{debnath}.  In an $n-$dimensional quasi Einstein manifold the Ricci tensor has precisely two distinct eigen values $a$ and $a+b$, where the multiplicity of $a$ is $n-1$ and $a+b$ is simple \cite{chacki1}. A proper $\eta-$Einstein contact metric manifold is a natural example of a quasi Einstein manifold (\cite{blair}, \cite{okumura}). Some of the physical and geometrical properties of quasi Einstein manifolds have been noticed in (\cite{29},  \cite{guha}, \cite{deghosh}, \cite{de5},  \cite{shaikh2},  \cite{dede}).

	The $k-$nullity distribution $N(k)$ of a Riemannian manifold $M_{n}$ is defined by
\begin{equation}
\label{1.2}
N(k): p{\longrightarrow} N_{p}(k)=\left\lbrace  Z\in{T_{p}M}:R(X,Y)Z =k\left[ g(Y,Z)X-g(X,Z)Y\right] \right\rbrace 
\end{equation}
for arbitrary vector fields $X$, $Y$ and $Z$, where $R$ represents the Riemannian curvature tensor  and $k$ is a smooth function on $M_{n}$ \cite{tanno}. The quasi Einstein manifold is called an $N(k)-$quasi Einstein manifold if the generator $\xi$ of the manifold $M_n$ belongs to $k-$nullity distribution  \cite{tripathi}.

	A space with constant curvature plays a central role in the development of differential geometry, Mathematical Physics and mechanics. Cartan (\cite{cartan1}, \cite{cartan2}) developed the idea of a locally symmetric Riemannian manifold in $1926$, which is a natural generalization of manifolds of constant curvature. The condition of local symmetry is equivalent to the fact that at every point $x \in M$, the local geodesic symmetry $F(x)$ is an isometry \cite{neill}. The idea of locally $\phi-$symmetric Sasakian manifold was introduced by Takahashi \cite{takahashi} in $1977$. Since then, the properties of such manifolds have been studied by several geometers on different spaces.

	The notion of a semi-symmetric linear connection on a differentiable manifold has been introduced by Friedmann and Schouten \cite{fs} in $1924$. Hayden \cite{hayden} in $1932$, introduced and studied the idea of semi-symmetric linear connection with torsion on a Riemannian manifold. After a long interval, Yano \cite{yano} started the systematic study of semi-symmetric metric connection on a Riemannian manifold in $1970$. In this connection, the properties of semi-symmetric metric connection are studied in (\cite{chaubey2}, \cite{sharfuddin}, \cite{imai}, \cite{rsmishra}, \cite{sunil1}, \cite{de2}, \cite{rnsingh}) and others. P. Zhao and H. Song \cite{zs} defined and studied a projective semi-symmetric connection on Riemannian manifold in $2001$. The properties of this connection has been studied by Zhao \cite{zhao}, Pal, Pandey and Singh \cite{pps} and others.

Motivated from the above studied, authors start the study of properties of Riemannian manifolds equipped with a projective semi-symmetric connection. The present paper is organized as follows:
After introductory section, we brief about the projective semi-symmetric connection in section $2$. It is prove that the curvature tensor with respect to projective semi-symmetric connection $\tilde{\nabla}$ coincide with the projective curvature tensor of the Levi-Civita connection $\nabla$. We also prove that the manifold $(M_n, g)$ endowed with $\tilde{\nabla}$ is a certain class of quasi Einstein manifold and the characteristic vector field $\xi$  belongs to $\lambda-$nullity distribution with respect to the connection $\tilde{\nabla}$ respectively. Section $3$ deals with the study of projective curvature tensor endowed with projective semi-symmetric connection $\tilde{\nabla}$ and prove that the projective curvature tensors with respect to connections $\tilde{\nabla}$ and $\nabla$ coincide. The properties of semi-symmetric Riemannian manifold admitting projective semi-symmetric connection $\tilde{\nabla}$ are given in section $4$. It is proved that the manifold is semi-symmetric for $\tilde{\nabla}$ if and only if it is flat. Section $5$ is the study about Riemannian manifold endowed with a projective semi-symmetric connection satisfying $\tilde{R}.\tilde{P}=0$ and prove some interesting results. In the last section we construct an example which support the existence of such connection $\tilde{\nabla}$ and verify some results.

\section{Projective semi-symmetric connection}
Let $M_n$ be an $n-$dimensional Riemannian manifold equipped with the Riemannian metric $g$ and $\nabla$ denotes the Levi-Civita connection on $(M_n, g)$ and satisfy 
\begin{equation}
\label{1}
\pi(X)=g(X, \xi) {\hspace{0.5cm}} and {\hspace{0.5cm}}g(\xi, \xi)=1
\end{equation}
for arbitrary vector field $X$, where $\pi$ is the first form associated with Riemannian metric $g$ and $\xi$ is an unit vector field of the manifold $(M_n,g)$. A linear connection $\tilde{\nabla}$ defined on $(M_n, g)$ is said to be semi-symmetric if the torsion tensor $\tilde{T}$ of the connection $\tilde{\nabla}$ defined as
$$\tilde{T}(X, Y)=\tilde{\nabla}_{X}Y-\tilde{\nabla}_{Y}X-[X, Y]$$ 
and satisfies 
$$\tilde{T}(X, Y)=\pi(Y)X-\pi(X)Y,$$
$i. e.$ $\tilde{T} \neq{0}$ for arbitrary vector fields $X$ and $Y$, otherwise it is symmetric. Also a linear connection $\tilde{\nabla}$ on $(M_n, g)$ is said to be metric if $\tilde{\nabla}g=0$, otherwise non-metric.  If the geodesics with respect to the linear connection $\tilde{\nabla}$ are consistent with those of Levi-Civita connection $\nabla$, then $\tilde{\nabla}$ is called the projective equivalent connection with $\nabla$. If we consider the linear connection $\tilde{\nabla}$ is semi-symmetric as well as projective equivalent, then it is called projective semi-symmetric connection \cite{zs}.  The tensorial relation between projective semi-symmetric and Levi-Civita connections on Riemannian manifold $(M_n, g)$ is given by
\begin{equation}
\label{2}
\tilde{\nabla}_{X}Y={\nabla}_{X}Y+\psi(Y)X+\psi(X)Y+\phi(Y)X-\phi(X)Y,
\end{equation}
for arbitrary vector fields $X$ and $Y$; where the $1-$form $\phi$ and $\psi$ are defined as 
\begin{equation}
\label{3}
\phi(X)=\frac{1}{2}\pi(X) {\hspace{0.5cm}}and {\hspace{0.5cm}} \psi(X)=\frac{n-1}{2(n+1)}\pi(X),
\end{equation}
for arbitrary vector fields $X$ and $Y$ \cite{zhao}. 
From (\ref{1}), (\ref{2}) and (\ref{3}), we have
\begin{equation}
\label{3a}
(\tilde{\nabla}_{X}g)(Y, Z)=\frac{1}{n+1}[2\pi(X)g(Y, Z)-n\pi(Y)g(X, Z)-n\pi(Z)g(X, Y)],
\end{equation}
for arbitrary vector fields $X$, $Y$ and $Z$. Thus the projective semi-symmetric connection $\tilde{\nabla}$ is non-metric. The properties of semi-symmetric non-metric connections have been noticed in (\cite{ac}, \cite{ch1}, \cite{ch2}, \cite{ch3}, \cite{ch4}, \cite{ch5}, \cite{ch6}) and many others.
It can be easily seen from \cite{zhao} that
\begin{equation}
\label{4}
\tilde{R}(X, Y)Z=R(X, Y)Z+\beta(X, Y)Z+\theta(X, Z)Y-\theta(Y, Z)X,
\end{equation}
where $\tilde{R}$ and $R$ denote the curvature tensors with respect to the connections $\tilde{\nabla}$ and $\nabla$ respectively and $\theta$, $\beta$ are $(0, 2)$ type tensors satisfying the following relations
\begin{equation}
\label{5}
\theta(X, Y)=(\nabla_{X}\phi)(Y)+(\nabla_{X}\psi)(Y)-\psi(X)\psi(Y)-\phi(X)\phi(Y)-\psi(X)\phi(Y)-\phi(X)\psi(Y),
\end{equation}
\begin{equation}
\label{6}
\beta(X, Y)=(\nabla_{X}\psi)(Y)-(\nabla_{Y}\psi)(X)+(\nabla_{X}\phi)(Y)-(\nabla_{Y}\phi)(X),
\end{equation}
for arbitrary vector fields $X$, $Y$ and $Z$. Let us suppose that the characteristic vector field $\xi$ is a parallel unit vector field with respect to the Levi-Civita connection, $i. e. , $ $\nabla{\xi}=0$ and $||\xi||=1$. This expression is equivalent to
\begin{equation}
\label{7}
(\nabla_{X}\pi)(Y)=\nabla_{X}\pi(Y)-\pi(\nabla_{X}Y)=0
\end{equation}
and therefore (\ref{2}) and (\ref{3}) assert that
\begin{equation}
\label{7a}
\tilde{\nabla}_{X}\xi=\frac{1}{n+1}\{nX-\pi(X)\xi\}.
\end{equation}
In consequence of (\ref{1}), (\ref{3}) and (\ref{7}), equations (\ref{5}) and (\ref{6}) become
\begin{equation}
\label{8}
\beta(X, Y)=0 {\hspace{0.25cm}} and  {\hspace{0.25cm}} \theta(X, Y)=\lambda \pi(X)\pi(Y),
\end{equation}
where $\lambda=-\frac{n^2}{(n+1)^2}$. It can be easily observe from (\ref{4}) and (\ref{8}) that
\begin{equation}
\label{9}
\tilde{R}(X, Y)Z=R(X, Y)Z+\lambda\{\pi(X)\pi(Z)Y-\pi(Y)\pi(Z)X\}.
\end{equation}
Contracting (\ref{9}) along $X$ and then using (\ref{1}), we have
\begin{equation}
\label{10}
\tilde{S}(Y,Z)=S(Y, Z)-\lambda (n-1) \pi(Y)\pi(Z)
\end{equation}
which gives
\begin{equation}
\label{11}
\tilde{r}=r-\lambda(n-1).
\end{equation}
Here $\tilde{S}$ and $S$; $\tilde{r}$ and $r$ denote the Ricci tensors and scalar curvatures corresponding to the connections $\tilde{\nabla}$ and $\nabla$ respectively.
\begin{theorem}
\label{newthm}
An $n-$dimensional Riemannian manifold $(M_n, g)$ equipped with a projective semi-symmetric connection $\tilde{\nabla}$  satisfying (\ref{7}) holds the following curvature conditions:
\newline
$(i) {\hspace{10.5pt}} '\tilde{R}(X, Y, Z, U)=- '\tilde{R}(Y, X, Z, U),$
\newline
$(ii) {\hspace{10.5pt}} '\tilde{R}(X, Y, Z, U) \neq - '\tilde{R}( X, Y, U, Z),$
\newline
$(iii) {\hspace{10.5pt}} '\tilde{R}(X, Y, Z, U) \neq '\tilde{R}(Z, U, X, Y),$
\newline
$(iv) {\hspace{10.5pt}} \tilde{R}(X, Y)Z+\tilde{R}(Y, Z)X+\tilde{R}(Z, X)Y=0,$
\newline
 $(v){\hspace{10.5pt}}(\tilde{\nabla}_{X}\tilde{R})(Y, Z)U+(\tilde{\nabla}_{Y}\tilde{R})(Z, X)U+(\tilde{\nabla}_{Z}\tilde{R})(X, Y)U $
\newline
$ {\hspace{25.5pt}}=2[\pi(X)R(Y, Z)U+\pi(Y)R(Z, X)U+\pi(Z)R(X, Y)U].$
\end{theorem}
\begin{proof}
From (\ref{9}), we have
\begin{equation}
\label{11_a}
 '\tilde{R}(X, Y, Z, U)='R(X, Y, Z, U)+\lambda\{\pi(X)\pi(Z)g(Y, U)-\pi(Y)\pi(Z)g(X, U)\},
\end{equation}
for all vector fields $X, Y, Z, U$ $\in T(M)$, where $'{\tilde{R}}(X, Y, Z, U)=g(\tilde{R}(X, Y) Z, U)$ and $'R(X, Y, Z, U)=g(R(X, Y) Z, U)$. By considering (\ref{11_a}) and curvature properties of $R$, we can easily verify the results (i), (ii) and (iii). With the help of (\ref{9}) and Bianchi's first identity, we get
\begin{equation}
\tilde{R}(X, Y)Z+\tilde{R}(Y, Z)X+\tilde{R}(Z, X)Y=0,
\end{equation}
which shows that the Riemann curvature tensor with respect to the projective semi-symmetric connection $\tilde{\nabla}$ satisfies the Bianchi's first identity. Covariant derivative of (\ref{9}) with respect to $\tilde{\nabla}$ gives
\begin{eqnarray}
\label{11a}
(\tilde{\nabla}_{X}\tilde{R})(Y, Z)U&=&(\tilde{\nabla}_{X}R)(Y, Z)U+\lambda\{\pi(U)(\tilde{\nabla}_{X}\pi)(Y)Z
+\pi(Y)(\tilde{\nabla}_{X}\pi)(U)Z\nonumber\\&&-\pi(U)(\tilde{\nabla}_{X}\pi)(Z)Y-\pi(Z)(\tilde{\nabla}_{X}\pi)(U)Y\}.
\end{eqnarray}
Also equations (\ref{1}), (\ref{2}), (\ref{3}), (\ref{3a}), (\ref{7}) and (\ref{7a}) yield
\begin{equation}
\label{11b}
 (\tilde{\nabla}_{X}\pi)(Y)=-\frac{n-1}{n+1}\pi(X)\pi(Y)
\end{equation}
and
\begin{eqnarray}
\label{11c}
(\tilde{\nabla}_{X}R)(Y, Z)U&=&({\nabla}_{X}R)(Y, Z)U+\frac{2}{n+1}\pi(X)R(Y, Z)U
-\frac{n}{n+1}\{\pi(Y)R(X, Z)U\nonumber\\&&+\pi(Z)R(Y, X)U+\pi(U)R(Y, Z)X\}.
\end{eqnarray}
In view of (\ref{11b}) and (\ref{11c}), (\ref{11a}) assumes the form
\begin{eqnarray}
\label{11d}
(\tilde{\nabla}_{X}\tilde{R})(Y, Z)U&=&({\nabla}_{X}R)(Y, Z)U+\frac{2}{n+1}\pi(X)R(Y, Z)U\nonumber\\&&
-\frac{n}{n+1}\{\pi(Y)R(X, Z)U+\pi(Z)R(Y, X)U+\pi(U)R(Y, Z)X\}\nonumber\\&&
-\frac{2\lambda(n-1)}{n+1}\{\pi(X)\pi(U)\pi(Y)Z-\pi(X)\pi(Z)\pi(U)Y\}.
\end{eqnarray}
The cyclic sum of (\ref{11d}) for the vector fields $X$, $Y$, $Z$ and use of Bianchi's second identity for $\nabla$ gives
\begin{eqnarray}
\label{11e}
&&(\tilde{\nabla}_{X}\tilde{R})(Y, Z)U+(\tilde{\nabla}_{Y}\tilde{R})(Z, X)U+(\tilde{\nabla}_{Z}\tilde{R})(X, Y)U\nonumber\\&&
=2[\pi(X)R(Y, Z)U+\pi(Y)R(Z, X)U+\pi(Z)R(X, Y)U].
\end{eqnarray}
This shows that a Riemannian manifold $(M_n, g)$ endowed with a projective semi-symmetric connection $\tilde{\nabla}$ satisfies the relation 
$$(\tilde{\nabla}_{X}\tilde{R})(Y, Z)U+(\tilde{\nabla}_{Y}\tilde{R})(Z, X)U+(\tilde{\nabla}_{Z}\tilde{R})(X, Y)U=0$$
if and only if 
$$\pi(X)R(Y, Z)U+\pi(Y)R(Z, X)U+\pi(Z)R(X, Y)U=0.$$
\end{proof}

	Let us suppose that the manifold $(M_n, g)$ is Ricci flat with respect to the projective semi-symmetric connection $\tilde{\nabla}$, $i. e.$, $\tilde{S}=0$ and therefore equation (\ref{10}) gives
\begin{equation}
\label{10a}
S(Y, Z)=\lambda (n-1)\pi(Y)\pi(Z),
\end{equation}
which shows that the manifold $(M_n, g)$, $(n>2)$, is a certain class of quasi Einstein manifold with the associated scalars $a=0$ and $b=\lambda (n-1)$.  
In consequence of (\ref{10a}), (\ref{9}) assumes the form
\begin{equation}
\label{10b}
\tilde{R}(X, Y)Z=P(X, Y)Z,
\end{equation}
where $P$ denotes the Weyl projective curvature tensor with respect to the Levi-Civita connection $\nabla$ and is given as
\begin{equation}
\label{10c}
P(X, Y)Z=R(X, Y)Z-\frac{1}{n-1}\{S(Y, Z)X-S(X, Z)Y\},
\end{equation}
for all vector fields $X, Y, Z$ $\in T(M)$. Thus we can conclude the results in the form of theorems as:
\begin{theorem}
\label{th1a}
Let $(M_n, g)$, $(n>2)$, be an $n-$dimensional Riemannian manifold equipped with a projective semi-symmetric connection $\tilde{\nabla}$ and the characteristic vector field $\xi$ of the manifold is a parallel unit vector field. If $(M_n, g)$ is Ricci flat with respect to the connection $\tilde{\nabla}$, then the projective curvature with respect to Levi-Civita connection $\nabla$ coincide with the curvature tensor of the connection $\tilde{\nabla}$.
\end{theorem}
\begin{theorem}
\label{thm1}
Let $(M_n, g)$, $(n>2)$, be an $n-$dimensional Riemannian manifold endowed with a projective semi-symmetric connection $\tilde{\nabla}$ and $\xi$ is a parallel unit vector field with respect to Levi-Civita connection. If $(M_n, g)$ is Ricci flat with respect to $\tilde{\nabla}$, then it is a certain class of quasi Einstein manifold.
\end{theorem}
By our assumption, the characteristic vector field $\xi$ is parallel unit vector field corresponding the Levi-Civita connection $\nabla$ and therefore by equation (\ref{7}) we can easily calculate that $R(X, Y)\xi=0$. After considering this fact and equation (\ref{1}), equation (\ref{9}) assumes the form
\begin{equation}
\label{12}
\tilde{R}(X, Y)\xi=\lambda\{\pi(X)Y-\pi(Y)X\}.
\end{equation}
This shows that the characteristic vector field $\xi$  belongs to the $\lambda-$nullity distribution with respect to the projective semi-symmetric connection $\tilde{\nabla}$. Thus we can state the following theorem:
\begin{theorem}
\label{thm2}
If $(M_n, g)$, $(n>2)$, be an $n-$dimensional Riemannian manifold admitting a projective semi-symmetric connection $\tilde{\nabla}$ and $\xi$ is a parallel unit vector field with respect to $\nabla$. Then the characteristic vector field of the manifold equipped with the projective semi-symmetric connection $\tilde{\nabla}$ belongs to $\lambda-$nullity distribution.
\end{theorem}
In view of  (\ref{1}), (\ref{7}),  (\ref{9}), (\ref{12}) and symmetric and skew-symmetric properties of curvature tensor,  we can state:
\begin{lemma}
\label{lem1}
If an $n-$dimensional Riemannian manifold $(M_n, g)$, $(n>2)$, admitting a projective semi-symmetric connection $\tilde{\nabla}$ and the characteristic vector field $\xi$ is a parallel unit vector field, then the following relations satisfy
\newline
$(i) {\hspace{0.25cm}}\tilde{R}(\xi, X)Y=\lambda\{\pi(Y)X-\pi(X)\pi(Y)\xi\}$,
\newline
$(ii){\hspace{0.25cm}} \tilde{R}(X, \xi)Y=\lambda \pi(Y)\{\pi(X)\xi-X\}$,
\newline
$(iii) {\hspace{0.25cm}}\pi(\tilde{R}(X, Y)Z)=0$,
for all vector fields $X, Y, Z$ $\in \chi(M)$. 
\end{lemma}
Proof is obvious by straight forward calculations.

 	In view of (\ref{1}), (\ref{7}) and (\ref{10}), we can compute that
\begin{equation}
\label{15}
(\tilde{\nabla}_{X}{\tilde{S}})(Y, Z)=(\nabla_{X}S)(Y, Z).
\end{equation}
Hence we can state the theorem:
\begin{theorem}
\label{thm4}
If $(M_n, g)$, $(n>2)$, be an $n-$dimensional Riemannian manifold equipped with a projective semi-symmetric connection $\tilde{\nabla}$ and $\xi$ is a parallel unit vector field. Then the manifold is Ricci-symmetric with respect to the projective semi-symmetric connection $\tilde{\nabla}$ if and only if it is Ricci-symmetric with respect to the Levi-Civita connection $\nabla$.
\end{theorem}
From equation (\ref{15}), we can also observe that
$$(\tilde{\nabla}_{X}{\tilde{S}})(Y, Z)-(\tilde{\nabla}_{Y}{\tilde{S}})(X, Z)=(\nabla_{X}S)(Y, Z)-(\nabla_{Y}S)(X, Z)$$ 
 and
\begin{eqnarray*}
&&(\tilde{\nabla}_{X}{\tilde{S}})(Y, Z)+(\tilde{\nabla}_{Y}{\tilde{S}})(Z, X)+(\tilde{\nabla}_{Z}{\tilde{S}})(X, Y)\nonumber\\&&
=(\nabla_{X}S)(Y, Z)+(\nabla_{Y}S)(Z, X)+(\nabla_{Z}S)(X, Y).
\end{eqnarray*}
 and hence the following lemma:
\begin{lemma}
\label{lem2}
Let $(M_n,g)$, $(n>2)$, be an $n-$dimensional Riemannian manifold endowed with a projective semi-symmetric connection $\tilde{\nabla}$ and $\xi$ is a parallel unit vector field. Then the Ricci tensor is of Codazzi type as well as cyclic parallel with respect to the projective semi-symmetric connection $\tilde{\nabla}$ if and only if it is Codazzi type as well cyclic parallel with respect to Levi-Civita connection $\nabla$ respectively.
\end{lemma}

\section{Projective curvature tensor equipped with projective semi-symmetric connection}
If  $\tilde{P}$ denotes the Weyl projective curvature tensor with respect to the connection $\tilde{\nabla}$, then
 \begin{equation}
\label{16}
^{\prime}{\tilde{P}}(X, Y, Z, U)=^{\prime}\tilde{R}(X, Y, Z, U)-\frac{1}{n-1}\{\tilde{S}(Y, Z)g(X, U)-\tilde{S}(X, Z)g(Y, U)\}
\end{equation}
holds for arbitrary vector fields $X$, $Y$, $Z$ and $U$, where $'{\tilde{P}}(X, Y, Z, U)=g(\tilde{P}(X, Y) Z, U)$.
In consequence of (\ref{9}) and (\ref{10}), above equation becomes
\begin{equation}
\label{17}
^{\prime}{\tilde{P}}(X, Y, Z, U)=^{\prime}{P}(X, Y, Z, U),
\end{equation}
where $P$ is the Weyl projective curvature tensor with respect to the Levi-Civita connection $\nabla$  given in (\ref{10c}) 
and $'P(X, Y, Z, U)=g(P(X, Y) Z, U)$. From the above discussions, we can conclude the result in the form of theorem as:
\begin{theorem}
\label{thm5}
If $(M_n, g)$, $(n>2)$, be an $n-$dimensional Riemannian manifold admitting a projective semi-symmetric connection $\tilde{\nabla}$ and $\xi$ is a parallel unit vector field, then the projective curvature tensors with respect to projective semi-symmetric and Levi-Civita connections coincide.
\end{theorem}
\begin{remark}
Zhao \cite{zhao} considered the special projective semi-symmetric connection and proved that the Weyl projective curvature tensors are invariant with respect to the special projective semi-symmetric and Levi-Civita connections. 
\end{remark}
From theorems \ref{th1a} and \ref{thm5}, we conclude the following:
\begin{theorem}
Let $(M_n, g)$, $(n>2)$, be an $n-$dimensional Riemannian manifold admitting a projective semi-symmetric connection $\tilde{\nabla}$ and $\xi$ is a parallel unit vector field. If the Ricci tensor with respect to projective semi-symmetric connection $\tilde{\nabla}$ is flat, then the curvature and Weyl projective curvature tensors with respect to  $\tilde{\nabla}$ and Weyl projective curvature tensor corresponding to Levi-Civita connection $\nabla$ coincide.
\end{theorem}

The Klein model of the Lobachevskii space and the central projection of a hemisphere on tangent space shows that the spaces of constant curvature are projectively flat and vice versa. Thus with the help of theorem \ref{thm5} and the above discussions, we can state:
\begin{theorem}
\label{thm6}
If $(M_n, g)$, $(n>2)$, be an $n-$dimensional Riemannian manifold admitting a projective semi-symmetric connection and $\xi$ is a parallel unit vector field, then $(M_n, g)$ is projectively flat with respect to the projective semi-symmetric connection if and only if it is of constant curvature.
\end{theorem}
Now, we consider that the Riemannian manifold is flat with respect to the projective semi-symmetric connection, $i. e., \tilde{R}=0$, then in consequence of theorem (\ref{th1a}) and equation (\ref{10b}), we obtain $P=0$. Therefore we can state the theorem:
\begin{theorem}
\label{thm7}
If an $n-$dimensional Riemannian manifold $(M_n, g)$, $(n>2)$, equipped with a projective semi-symmetric connection and a parallel unit vector field is flat with respect to projectively semi-symmetric connection $\tilde{\nabla}$, then it is a manifold of constant curvature although converse part is also true.
\end{theorem}
\begin{remark}
The idea of constant curvature is playing a central role in the theory of relativity and cosmology. The simplest cosmological model can be constructed by assuming that the universe is isotropic and homogeneous. This is known as cosmological principle. When we translated this principle to Riemannian geometry, professes that the three dimensional position space is a space of maximal symmetry \cite{stephani}, $i. e.$, a space of constant curvature whose curvature depends upon time. The cosmological solutions of Einstein equations which contain a three dimensional space like surfaces of a constant curvature are the Robertson-Walker metric, while four dimensional space of constant curvature is the de Sitler model of the universe [\cite{narlikar}, \cite{stephani}]. 
\end{remark}

\section{semi-symmetric Riemannian manifold admitting projective semi-symmetric connection }
A Riemannian manifold $(M_n, g)$ is said to be semi-symmetric \cite{szabo1}, \cite{szabo2} with respect to the Levi-Civita connection $\nabla$ if its non-flat curvature tensor $R$ satisfies the condition $R. R=0$ . Analogous to this, we can define:
\begin{definition}
A non-flat Riemannian manifold $(M_n, g)$, $(n>2)$, is said to be semi-symmetric with respect to the projective semi-symmetric connection $\tilde{\nabla}$ if $\tilde{R}.\tilde{R}=0$. 
\end{definition}

It is obvious that
\begin{eqnarray*}
\label{19}
(\tilde{R}(\xi, X).\tilde{R})(Y, Z)U&=&\tilde{R}(\xi, X)\tilde{R}(Y, Z)U-\tilde{R}(\tilde{R}(\xi, X)Y, Z)U\nonumber\\&&
-\tilde{R}(Y, \tilde{R}(\xi, X)Z)U-\tilde{R}(Y, Z)\tilde{R}(\xi, X)U.
\end{eqnarray*}
In view of (\ref{1}), (\ref{12}) and  lemma (\ref{lem1}), last equation becomes
\begin{eqnarray}
\label{20}
(\tilde{R}(\xi, X).\tilde{R})(Y, Z)U&=&-\lambda\{\pi(Y)\tilde{R}(X, Z)U+\pi(Z)\tilde{R}(Y, X)U+\pi(U)\tilde{R}(Y, Z)X\}\nonumber\\&&
+2\lambda^2\{\pi(Y)Z-\pi(Z)Y\}\pi(X)\pi(U).
\end{eqnarray}
Let us suppose that $\tilde{R}.\tilde{R}=0$, then equation (\ref{20}) reflects that either $\lambda=0$ or $\pi(Y)\tilde{R}(X, Z)U+\pi(Z)\tilde{R}(Y, X)U+\pi(U)\tilde{R}(Y, Z)X=2\lambda\{\pi(Y)Z-\pi(Z)Y\}\pi(X)\pi(U)$. Since $\lambda \neq {0}$, therefore
$$\pi(Y)\tilde{R}(X, Z)U+\pi(Z)\tilde{R}(Y, X)U+\pi(U)\tilde{R}(Y, Z)X=2\lambda\{\pi(Y)Z-\pi(Z)Y\}\pi(X)\pi(U).$$
Putting $U=\xi$ in last expression and then using (\ref{1}) and (\ref{12}), we find
\begin{equation}
\label{21}
\tilde{R}(Y, Z)X=\lambda \pi(X)\{\pi(Y)Z-\pi(Z)Y\},
\end{equation} 
which is equivalent to $R(Y, Z)X=0$. Thus the manifold $(M_n, g)$, $(n>2)$, equipped with a projective semi-symmetric connection $\tilde{\nabla}$ satisfying $\tilde{R}.\tilde{R}=0$ is flat for Levi-Civita connection $\nabla$. Conversely, if the manifold is flat, $i. e.$, $R=0$, then equation (\ref{9}) assumes the form (\ref{21}).  With the help of (\ref{1}), (\ref{12}) and (\ref{21}), (\ref{20}) shows that $\tilde{R}.\tilde{R}=0$.
Hence we state the above result in the form of theorem as:
\begin{theorem}
\label{thm8}
Let $(M_n, g)$, $(n>2)$, be an $n-$dimensional Riemannian manifold admitting a projective semi-symmetric connection $\tilde{\nabla}$ satisfying (\ref{7}). Then the necessary and sufficient condition for a manifold $(M_n, g)$ to be semi-symmetric with respect to the connection $\tilde{\nabla}$ is that the manifold is flat with respect to Levi-Civita connection $\nabla$.
\end{theorem}
In consequence of (\ref{9}), (\ref{11d}),  $\lambda=-\frac{n^2}{(n+1)^2}$ and theorem {\ref{thm8}}, we find
\begin{equation*}
(\tilde{\nabla}_{X}\tilde{R})(Y, Z)U=\rho(X)\tilde{R}(Y, Z)U,
\end{equation*}
where $\rho(X)= -\frac{2(n-1)}{n+1}\pi(X)$. Thus we can state:
\begin{corollary}
\label{cor_1}
Let $(M_n, g)$, $(n>2)$, be an $n-$dimensional Riemannian manifold admitting a projective semi-symmetric connection $\tilde{\nabla}$ and satisfies (\ref{7}). If $(M_n, g)$ is semi-symmetric with respect to the projective semi-symmetric connection $\tilde{\nabla}$, then it is recurrent.
\end{corollary}
From theorem {\ref{thm8}}, we can also conclude the following corollary: 
\begin{corollary}
\label{cor_2}
A semi-symmetric Riemannian manifold $(M_n, g)$, $(n>2)$, endowed with a projective semi-symmetric connection  $\tilde{\nabla}$ is projectively, conformally, concircularly, conharmonicaly, quasi-conformally and $m-$projectively flat for Levi-Civita connection $\nabla$. 
\end{corollary}
In view of theorem {\ref{thm5}} and corollary {\ref{cor_2}}, we have
\begin{corollary}
Every semi-symmetric Riemannian manifold $(M_n, g)$, $(n>2)$, equipped with a projective semi-symmetric connection  $\tilde{\nabla}$ is projectively flat for  $\tilde{\nabla}$.
\end{corollary}

\section{Riemannian manifold equipped with a projective semi-symmetric connection satisfying $\tilde{R}.\tilde{P}=0$}
We have
\begin{eqnarray*}
\label{5.1}
(\tilde{R}(X, Y).\tilde{P})(Z, U)V&=&\tilde{R}(X, Y)\tilde{P}(Z, U)V-\tilde{P}(\tilde{R}(X, Y)Z, U)V\nonumber\\&&
-\tilde{P}(Z, \tilde{R}(X, Y)U)V-\tilde{P}(Z, U)\tilde{R}(X, Y)V.
\end{eqnarray*}
Replacing $X$ by $\xi$ in above equation and then using equation (\ref{1}) and lemma (\ref{lem1}), we have 
\begin{eqnarray}
\label{5.1}
(\tilde{R}(\xi, Y).\tilde{P})(Z, U)V&=&\lambda[\pi(\tilde{P}(Z, U)V)Y-\pi(\tilde{P}(Z, U)V)\pi(Y)\xi-\pi(Z)\tilde{P}(Y, U) V\nonumber\\&&
-\pi(U)\tilde{P}(Z, Y)V-\pi(V)\tilde{P}(Z, U)Y+\pi(Y)\pi(Z)\tilde{P}(\xi, U)V\nonumber\\&&
+\pi(U)\pi(Y)\tilde{P}(Z, \xi)V+\pi(Y)\pi(V)\tilde{P}(Z, U)\xi].
\end{eqnarray}
Let us suppose that $\tilde{R}.\tilde{P}=0$, then we have from (\ref{5.1}) ($\lambda \neq {0}$)
\begin{eqnarray*}
&&\pi(Z)\tilde{P}(Y, U) V+\pi(U)\tilde{P}(Z, Y)V+\pi(V)\tilde{P}(Z, U)Y\nonumber\\&&
=\pi(\tilde{P}(Z, U)V)Y-\pi(\tilde{P}(Z, U)V)\pi(Y)\xi+\pi(Y)\pi(Z)\tilde{P}(\xi, U)V\nonumber\\&&
+\pi(U)\pi(Y)\tilde{P}(Z, \xi)V+\pi(Y)\pi(V)\tilde{P}(Z, U)\xi.
\end{eqnarray*}
Setting $Z$ by $\xi$ in the above equation and using (\ref{1}), we get
\begin{eqnarray}
\label{5.2}
&&\tilde{P}(Y, U) V+\pi(U)\tilde{P}(\xi, Y)V+\pi(V)\tilde{P}(\xi, U)Y\nonumber\\&&
=\pi(\tilde{P}(\xi, U)V)Y-\pi(\tilde{P}(\xi, U)V)\pi(Y)\xi+\pi(Y)\tilde{P}(\xi, U)V\nonumber\\&&
+\pi(U)\pi(Y)\tilde{P}(\xi, \xi)V+\pi(Y)\pi(V)\tilde{P}(\xi, U)\xi.
\end{eqnarray}
 In consequence of (\ref{1}), (\ref{10}), lemma {\ref{lem1}} and (\ref{16}), we find that
\begin{eqnarray}
\label{5.3}
&&(i) {\hspace{10.5pt}}\tilde{P}(\xi, X)Y=\frac{1}{n-1}\{S(\xi, Y)X-S(X, Y)\xi\},
\nonumber\\&&
(ii){\hspace{10.5pt}} \pi(\tilde{P}(X, Y)Z)=\frac{1}{n-1}\{\pi(Y)S(X, Z)-\pi(X)S(Y, Z)\}.
\end{eqnarray}
In view of (\ref{5.3}), (\ref{5.2}) turns into the form
\begin{eqnarray}
\label{5.4}
\tilde{P}(Y, U) V&=&\frac{1}{n-1}\{\pi(U)S(V, Y)\xi-\pi(V)S(\xi, Y)U+\pi(V)S(U, Y)\xi
\nonumber\\&&
-\pi(U)\pi(Y)S(\xi, V)\xi-S(U, V)Y+\pi(Y)S(\xi, V)U
\nonumber\\&&
+\pi(Y)\pi(V)S(\xi, \xi)U-\pi(Y)\pi(V)S(U, \xi)\xi\}.
\end{eqnarray}
which gives
\begin{eqnarray}
\label{5.4a}
\pi(\tilde{P}(Y, U) V)&=&\frac{1}{n-1}\{\pi(U)S(V, Y)-\pi(V)S(\xi, Y)\pi(U)+\pi(V)S(U, Y)
\nonumber\\&&
-\pi(U)\pi(Y)S(\xi, V)-S(U, V)\pi(Y)+\pi(Y)S(\xi, V)\pi(U)
\nonumber\\&&
+\pi(Y)\pi(V)S(\xi, \xi)\pi(U)-\pi(Y)\pi(V)S(U, \xi)\}.
\end{eqnarray}
Using (\ref{5.3}) $(ii)$ in (\ref{5.4a}), we obtain
\begin{eqnarray*}
\label{5.4b}
\pi(V)S(U, Y)&=&\pi(U)\pi(V)S(\xi, Y)+\pi(U)\pi(Y)S(\xi, V)-\pi(U)\pi(Y)S(\xi, V)\nonumber\\&&
-\pi(Y)\pi(V)\pi(U)S(\xi, \xi)+\pi(Y)\pi(V)S(U, \xi).
\end{eqnarray*}
Setting $V=\xi$ in above equation and use of (\ref{1}) gives
\begin{equation}
\label{5.4c}
S(U, Y)=\pi(U)S(\xi, Y)+\pi(Y)S(U, \xi)-\pi(U)\pi(Y)S(\xi, \xi).
\end{equation}
Using (\ref{1}), (\ref{5.4c}) in (\ref{5.4}), we find
\begin{eqnarray}
\label{5.4d}
'{\tilde{P}}(Y, U, V, X)&=&\frac{1}{n-1}\{2\pi(U)\pi(X)\pi(V)S(\xi, Y)-2\pi(X)\pi(U)\pi(V)\pi(Y)S(\xi, \xi)\nonumber\\&&
-\pi(V)S(\xi, Y)g(U, X)-\pi(U)S(\xi, V)g(Y, X)\nonumber\\&&
-\pi(V)S(U, \xi)g(X, Y)+\pi(U)\pi(V)g(X, Y)S(\xi, \xi)\nonumber\\&&
+\pi(Y)S(\xi, V)g(U, X)+\pi(Y)\pi(V)S(\xi, \xi)g(U, X)\}.
\end{eqnarray}
Let $\left\{\left\{e_{i}\right\},  i=1, 2,...,n\right\}$, be an orthonormal basis of the tangent space at any point of the manifold $(M_n, g)$. Then putting $X=Y=e_{i}$ in (\ref{5.4d}) and taking summation over $i$, $1 \leq i \leq n$, we get 
\begin{equation}
\label{5.5}
\pi(U)S(\xi, V)=\frac{n+1}{n-1}\{\pi(U)\pi(V)S(\xi, \xi)-\pi(V)S(U, \xi)\}.
\end{equation}
after considering equations (\ref{1}), (\ref{16}) and $\sum_{i=1}^{n} {'{\tilde{P}}(e_i, U, V, e_i)}=0$. Replacing the vector field $U$ with $\xi$ in (\ref{5.5}) and use of (\ref{1}) gives
\begin{equation*}
\label{5.6}
S(\xi, V)=0
\end{equation*}
 and therefore equation (\ref{5.4c}) gives
$$S(U, Y)=0.$$
for arbitrary vector fields $U$ and $Y$. Above equation shows that the manifold $(M_n, g)$, $(n>2)$, is Ricci flat with respect to the Levi-Civita connection. Thus we can state:
\begin{theorem}
\label{thm_last}
A Riemannian manifold $(M_n, g)$, $(n>2)$, equipped with a projective semi-symmetric connection $\tilde{\nabla}$  satisfying $\tilde{R}.\tilde{P}=0$  is Ricci flat.
\end{theorem}
By the use of (\ref{1}), (\ref{9}), (\ref{10}), (\ref{16}), theorems \ref{thm5} and \ref{thm_last}), we conclude that
\begin{equation*}
\tilde{P}(Z, U)V=R(Z, U)V=P(Z, U)V.
\end{equation*}
\begin{corollary}
If a Riemannian manifold $(M_n, g)$ admits a projective semi-symmetric connection $\tilde{\nabla}$ and satisfies the condition $\tilde{R}.\tilde{P}=0$, then the curvature tensor  $R$ and projective curvature tensor $P$ for $\nabla$ coincides with projective curvature tensor of $\tilde{\nabla}$.
\end{corollary}
From theorems \ref{thm8} and \ref{thm_last}, we conclude the following corollary:
\begin{corollary} 
Let $(M_n, g)$, $(n>2)$, be an $n-$dimensional Riemannian manifold equipped with a projective semi-symmetric connection $\tilde{\nabla}$ and $\xi$ is a parallel unit vector field. Then every semi-symmetric Riemannian manifold with respect to $\tilde{\nabla}$ satisfies $\tilde{R}.\tilde{P}=0$.
\end{corollary}

\section{Example}
 P. Alegre, D. E. Blair and A. Carriazo \cite{alegre} introduced the idea of generalized Sasakian space form and they constructed many examples by using some different geometric techniques such as Riemannian submersions, warped products or conformal and related transformations in $2004$,. A Riemannian manifold $M_n$ of dimension $n$ equipped with a tensor field $\phi$ of type $(1, 1)$, a structure vector field $\xi$ and a covariant vector field $\eta$ associated with the Riemannian metric $g$ satisfies the relations
\begin{equation} 
\label{*1}
\phi^2(X)=-X+\eta(X)\xi, {\hspace{0.25cm}} \eta(\xi)=1, {\hspace{0.25cm}} g(X, \xi)=\eta(X), {\hspace{0.25cm}} and   {\hspace{0.5cm}} \phi{\xi}=0
\end{equation}
and
\begin{equation*} 
\label{*1a}
g(X, Y)=g(\phi{X}, \phi{Y})+\eta(X)\eta(Y),
\end{equation*}
for arbitrary vector fields $X$ and $Y$, is called an almost contact metric manifold $(M_n, \phi, \xi, \eta, g)$ \cite{blair}. An almost contact metric manifold $(M_n, \phi, \xi, \eta, g)$ is cosymplectic \cite{blair} if $\nabla{\phi}=0$, which implies the following expressions
\begin{equation*}
\nabla{\xi}=0, {\hspace{0.25cm}} \nabla{\eta}=0  {\hspace{0.25cm}}and {\hspace{0.25cm}}R(X, Y)\xi=0.
\end{equation*}
An almost contact metric manifold $(M_n, \phi, \xi, \eta, g)$ is said to be a generalized Sasakian space form \cite{alegre} if the Riemannian curvature tensor $R$ satisfies the tensorial relation
\begin{eqnarray}
\label{*2}
&&R(X, Y)Z=f_{1}\{g(Y, Z)X-g(X, Z)Y\}\nonumber\\&&
+f_{2}\{g(X, \phi{Z})\phi{Y}-g(Y, \phi{Z})\phi{X}+2g(X, \phi{Y})\phi{Z}\}\nonumber\\&&
+f_{3}\{\eta(X)\eta(Z)Y-\eta(Y)\eta(Z)X+g(X, Z)\eta(Y)\xi-g(Y, Z)\eta(X)\xi\}
\end{eqnarray}
for arbitrary vector fields $X$, $Y$ and $Z$, where $f_{1}$, $f_2$ and $f_3$ are smooth functions on $M_n$. First author proved that the generalized Sasakian space form is a certain class of quasi Einstein manifold \cite{chaubey1add}. Here we suppose that the manifold $M_n$ is cosymplectic as well as generalized Sasakian space form and $f_1 = f_3 \neq {0}$. Replacing $Z$ with the structure vector field $\xi$  in (\ref{*2}) and then using equation (\ref{*1}),  we get
\begin{equation}
\label{*3}
R(X, Y)\xi=0.
\end{equation}
If we define the projective semi-symmetric connection $\tilde{\nabla}$ on $M_n$  as
\begin{equation}
\label{*4}
\tilde{\nabla}_{X}Y={\nabla}_{X}Y+\frac{n-1}{2(n+1)}\{\eta(Y)X+\eta(X)Y\}+\frac{1}{2}\{\eta(Y)X-\eta(X)Y\},
\end{equation}
and torsion tensor of the connection $\tilde{\nabla}$
\begin{equation*}
\tilde{T}(X, Y)=\eta(Y)X-\eta(X)Y
\end{equation*}
for arbitrary vector fields $X$ and $Y$, then with the help of equations (\ref{9}), (\ref{*3}) and (\ref{*4}), we obtain the expression (\ref{12}) and hence the theorem \ref{thm1} verified.


\end{document}